\documentclass{amsart}
\usepackage{amsmath,amsxtra,amssymb}
\usepackage[all]{xy}

\newcommand{\hz}{\vspace{0.5cm}}

\renewcommand{\qed}{\hspace*{\fill}$\Box$\hz\pagebreak[1]}

\newcommand{\M}{{\mathcal M}}

\newtheorem{theorem}{Theorem}[section]
\newtheorem{cor}[theorem]{Corollary}
\newtheorem{prop}[theorem]{Proposition}
\newtheorem{rem}[theorem]{Remark}
\newtheorem{lemma}[theorem]{Lemma}

\newtheorem{exam}[theorem]{Example}

\begin{document}
\title{The Effros-Ruan conjecture for bilinear forms on C$^*$-algebras}
\author{Uffe Haagerup$^{(1)}$ and Magdalena Musat$^{(2)}$}
\address{$^{(1)}$ Department of Mathematics and
Computer Science, University of Southern Denmark, Campusvej 55, 5230
Odense M, Denmark.\\
$^{(2)}$ Department of Mathematical Sciences, University of Memphis,
373 Dunn Hall, Memphis, TN, 38152, USA.}
\email{$^{(1)}$haagerup@imada.sdu.dk\\$^{(2)}$mmusat@memphis.edu}.

%\author{Magdalena Musat}
%\address{Department of Mathematics, 0112\\
%University of California, San Diego\\
%La Jolla, CA 92093-0112}
%\email{mmusat@math.ucsd.edu}
\date{}

\footnotetext {$^{(1)}$ Partially supported by the Danish Natural
Science Research Council.\\ \hspace*{0.48cm} $^{(2)}$ Partially
supported by the National Science Foundation, DMS-0703869.}


\maketitle

\begin{abstract}
In 1991 Effros and Ruan conjectured that a certain Grothendieck-type
inequality for a bilinear form on C$^*$-algebras holds if (and only
if) the bilinear form is jointly completely bounded. In 2002 Pisier
and Shlyakhtenko proved that this inequality holds in the more
general setting of operator spaces, provided that the operator
spaces in question are exact. Moreover, they proved that the
conjecture of Effros and Ruan holds for pairs of C$^*$-algebras, of
which at least one is exact. In this paper we prove that the
Effros-Ruan conjecture holds for general C$^*$-algebras, with
constant one. More precisely, we show that for every jointly
completely bounded (for short, j.c.b.) bilinear form on a pair of
C$^*$-algebras $A$ and $B$\,, there exist states $f_1$\,, $f_2$ on
$A$ and $g_1$\,, $g_2$ on $B$ such that for all $a\in A$ and $b\in
B$\,,
\begin{equation*}
|u(a, b)|\leq
\|u\|_{jcb}(f_1(aa^*)^{1/2}g_1(b^*b)^{1/2}+f_2(a^*a)^{1/2}g_2(bb^*)^{1/2})\,.
\end{equation*}
While the approach by Pisier and Shlyakhtenko relies on free
probability techniques, our proof uses more classical operator
algebra theory, namely, Tomita-Takesaki theory and special
properties of the Powers factors of type III$_\lambda$\,, $0<
\lambda< 1$\,.
\end{abstract}

%\maketitle
\section{Introduction}
\setcounter{equation}{0}

In 1956 Grothendieck published the celebrated "R\'esum\'e de la
th\'eorie m\'etrique des produits tensoriels topologiques",
containing a general theory of tensor norms on tensor products of
Banach spaces, describing several operations to generate new norms
from known ones, and studying the duality theory between these
norms. Since 1968 it has had considerable influence on the
development of Banach space theory (see e.g., \cite{LT})\,. The
highlight of the paper \cite{Gro}, now referred to as the
"R\'esum\'e" is a result that Grothendieck called "The fundamental
theorem on the metric theory of tensor products". Grothendieck's
theorem asserts that given compact spaces $K_1$ and $K_2$ and a
bounded bilinear form $u:C(K_1)\times C(K_2)\rightarrow \mathbb{K}$
(where $\mathbb{K}=\mathbb{R}$ or $\mathbb{C}$)\,, then there exist
probability measures $\mu_1$ and $\mu_2$ on $K_1$ and $K_2$\,,
respectively, such that
\[ |u(f, g)|\leq K_G^{\mathbb{K}} \|u\|\left(\int_{K_1} |f(t)|^2
\,d\mu_1(t)\right)^{1/2} \left(\int_{K_2}
|g(t)|^2\,d\mu_2(t)\right)^{1/2}\,,
\]
for all $f\in C(K_1)$ and $g\in C(K_2)$\,, where $K_G^{\mathbb{K}}$
is a universal constant.

The non-commutative version of Grothendieck's inequality
(conjectured in the "R\'esum\'e") was first proved by Pisier under
some approximability assumption (cf.\ \cite{Pi7})\,, and obtained in
full generality in \cite{Ha6}. The theorem asserts that given
C$^*$-algebras $A$ and $B$ and a bounded bilinear form $u:A\times
B\rightarrow \mathbb{C}$\,, then there exist states $f_1\,, f_2$ on
$A$ and states $g_1\,, g_2$ on $B$ such that for all $a\in A$ and
$b\in B$\,,
\begin{equation*}
|u(a, b)|\leq \|u\|(f_1(a^*a)+f_2(aa^*))^{1/2}(g_1(b^*b)+
g_2(bb^*))^{1/2}\,.
\end{equation*}
As a corollary, it was shown in \cite{Ha6} that given C$^*$-algebras
$A$ and $B$\,, then any bounded linear operator $T: A\rightarrow
B^*$ admits a factorization $T=SR$ through a Hilbert space $H$\,,
where $A \overset{R}{\longrightarrow} H \overset{S}{\longrightarrow}
B^*$\,, and
\begin{equation*}
\|R\|\|S\|\leq 2\|T\|\,.
\end{equation*}
Let $E\subseteq A$ and $F\subseteq B$ be operator spaces sitting in
C$^*$-algebras $A$ and $B$\,, and let $u: E\times F\rightarrow
\mathbb{C}$ be a bounded bilinear form. Then, there exists a unique
bounded linear operator $\widetilde{u}:E\rightarrow F^*$ such that
\begin{equation}\label{eq555555567}
u(a\,, b):=\langle \widetilde{u}(a)\,, b\rangle\,, \quad a\in E\,,
b\in F\,,
\end{equation}
where $\langle \,\cdot\,, \cdot\, \rangle$ denotes the duality
bracket between $F$ and $F^*$\,. The map $u$ is called {\em jointly
completely bounded} (for short, j.c.b.) if the associated map
$\widetilde{u}:E\rightarrow F^*$ is completely bounded, in which
case we set
\begin{equation}\label{eq77777775675}
\|u\|_{\text{jcb}}:=\|\widetilde{u}\|_{\text{cb}}\,. \end{equation}
(Otherwise, we set $\|u\|_{\text{jcb}}=\infty$\,.) It is easily
checked that
\begin{equation}\label{eq66655566543}
\|u\|_{\text{jcb}}=\sup\limits_{n\in \mathbb{N}} \|u_n\|\,,
\end{equation}
where for every $n\geq 1$\,, the map $u_n: M_n(E)\otimes
M_n(F)\rightarrow M_n(\mathbb{C})\otimes M_n(\mathbb{C})$ is given
by
\begin{equation*}
u_n\left(\sum_{i=1}^k a_i\otimes c_i\,, \sum_{j=1}^l b_j\otimes d_j
\right)=\sum_{i=1}^k \sum_{j=1}^l u(a_i, b_j)c_i\otimes d_j\,,
\end{equation*}
for all finite sequences $\{a_i\}_{1\leq i\leq k}$ in $E$\,,
$\{b_j\}_{1\leq j\leq l}$ in $F$\,, $\{c_i\}_{1\leq i\leq k}$ and
$\{d_j\}_{1\leq j\leq l}$ in $M_n(\mathbb{C})$\,, $k\,, l\in
\mathbb{N}\,.$ Moreover, $\|u\|_{\text{jcb}}$ is the smallest
constant $\kappa_1$ for which, given arbitrary C$^*$-algebras $C$
and $D$ and finite sequences $\{a_i\}_{1\leq i\leq k}$ in $E$\,,
$\{b_j\}_{1\leq j\leq l}$ in $F$\,, $\{c_i\}_{1\leq i\leq k}$ in $C$
and $\{d_j\}_{1\leq j\leq l}$ in $D$\,, where $k\,, l\in
\mathbb{N}$\,, the following inequality holds
\begin{equation}\label{eq3333222333220}
\left\|\sum_{i=1}^k \sum_{j=1}^l u(a_i\,, b_j) c_i\otimes
d_j\right\|_{C\otimes_{\text{min}} D}\leq \kappa_1
\left\|\sum_{i=1}^k a_i\otimes c_i\right\|_{E\otimes_{\text{min}} C}
\left\|\sum_{j=1}^l b_j\otimes d_j\right\|_{F\otimes_{\text{min}}
D}\,.
\end{equation}
For a reference, see the discussion following Definition 1.1 in
\cite{PiSh}

It was conjectured by Effros and Ruan in 1991 (cf.\ \cite{ER2} and
\cite{PiSh}, Conjecture 0.1) that if $A$ and $B$ are C$^*$-algebras
and $u:A\times B\rightarrow \mathbb{C}$ is a jointly completely
bounded bilinear form, then there exist states $f_1$\,, $f_2$ on $A$
and states $g_1$\,, $g_2$ on $B$ such that for all $a\in A$ and
$b\in B$\,,
\begin{equation}\label{eq2221112221113}
|u(a, b)|\leq
K\|u\|_{jcb}(f_1(aa^*)^{1/2}g_1(b^*b)^{1/2}+f_2(a^*a)^{1/2}g_2(bb^*)^{1/2})\,,
\end{equation}
where $K$ is a universal constant.

In \cite{PiSh} Pisier and Shlyakhtenko proved an operator space
version of (\ref{eq2221112221113}), namely, if $E\subseteq A$ and
$F\subseteq B$ are exact operator spaces with exactness constants
$\text{ex}(E)$ and $\text{ex}(F)$\,, respectively, and $u: E\times
F\rightarrow \mathbb{C}$ is a j.c.b. bilinear form, then there exist
states $f_1$\,, $f_2$ on $A$ and states $g_1$\,, $g_2$ on $B$ such
that for all $a\in E$ and $b\in F$\,,
\begin{equation*}
|u(a, b)|\leq
2^{3/2}\text{ex}(E)\text{ex}(F)\|u\|_{jcb}(f_1(aa^*)^{1/2}g_1(b^*b)^{1/2}+f_2(a^*a)^{1/2}g_2(bb^*)^{1/2})\,.
\end{equation*}
Moreover, by the same methods they were able to prove the
Effros-Ruan conjecture for C$^*$-algebras with constant
$K=2^{3/2}$\,, provided that at least one of the C$^*$-algebras
$A$\,, $B$ is exact (cf.\ \cite{PiSh}, Theorem 0.5).

The main result of this paper is that the Effros-Ruan conjecture is
true. Moreover, it holds with constant $K=1$\,, that is,

\begin{theorem}\label{theffruan}
Let $A$ and $B$ be C$^*$-algebras and $u:A\times B\rightarrow
\mathbb{C}$ a jointly completely bounded bilinear form. Then there
exist states $f_1$\,, $f_2$ on $A$ and states $g_1$\,, $g_2$ on $B$
such that for all $a\in A$ and $b\in B$\,,
\begin{equation*}
|u(a, b)|\leq
\|u\|_{jcb}(f_1(aa^*)^{1/2}g_1(b^*b)^{1/2}+f_2(a^*a)^{1/2}g_2(bb^*)^{1/2})\,.
\end{equation*}
\end{theorem}

It follows from Theorem \ref{theffruan} that every completely
bounded linear map $T:A\rightarrow B^*$ from a C$^*$-algebra $A$ to
the dual $B^*$ of a C$^*$-algebra $B$ has a factorization $T=vw$
through $H_r\oplus K_c$ (the direct sum of a row Hilbert space and a
column Hilbert space), such that
\begin{equation*}
\|v\|_{\text{cb}}\|w\|_{\text{cb}}\leq 2\|T\|_{\text{cb}}\,.
\end{equation*}
(See Proposition \ref{prop565656445} of this paper.) Theorem
\ref{theffruan} also settles in the affirmative a related conjecture
by Blecher (cf.\ \cite{Bl}; see also \cite{PiSh}, Conjecture 0.2).
For details, see Remark \ref{rem56743233332} of this paper.

Furthermore, thanks to Theorem \ref{theffruan} we can strengthen a
number of results from \cite{PiSh}\,, cf.\ Corollaries 3.7 through
3.10 in this paper. For instance, it follows that if an operator
space $E$ and its dual $E^*$ both embed in noncommutative
$L_1$-spaces, then $E$ is completely isomorphic to a quotient of a
subspace of $H_r\oplus K_c$\,, for some Hilbert spaces $H$ and
$K$\,.

It also follows from Theorem \ref{theffruan} that if $u:A\times
B\rightarrow \mathbb{C}$ is a j.c.b. bilinear form on C$^*$-algebras
$A$ and $B$\,, then the inequality (\ref{eq3333222333220}) holds, as
well, when the $C\otimes_{\text{min}} D$-norm on the left-hand side
is replaced by the $C\otimes_{\text{max}} D$-norm (with constant
$2\|u\|_{\text{jcb}}$ instead of $\|u\|_{\text{jcb}}$)\,, cf.\
Proposition \ref{prop789}. Moreover, we show that for bilinear forms
$u$ on operator spaces $E\subseteq A$ and $F\subseteq B$ sitting in
C$^*$-algebras $A$ and $B$\,, the above mentioned variant of
(\ref{eq3333222333220})\,, namely the inequality
\begin{equation*}
\left\|\sum\limits_{i=1}^m \sum\limits_{j=1}^n u(a_i, b_j)c_i\otimes
d_j \right\|_{C\otimes_{\text{max}} D}\leq \kappa_4
\left\|\sum\limits_{i=1}^m a_i\otimes
c_i\right\|_{E\otimes_{\text{min}} C}\left\|\sum\limits_{j=1}^n
b_j\otimes d_j\right\|_{F\otimes_{\text{min}} D}
\end{equation*}
(where $C$ and $D$ are arbitrary C$^*$-algebras) characterizes those
j.c.b. bilinear forms that satisfy an Effros-Ruan type inequality\,.
That is, there exists a constant $\kappa_2\geq 0$ and states $f_1\,,
f_2$ on $A$ and states $g_1\,, g_2$ on $B$ such that, for all $a\in
E$ and $b\in F$\,,
\begin{equation*}
 |u(a, b)|\leq \kappa_2
(f_1(aa^*)^{1/2}g_1(b^*b)^{1/2}+f_2(a^*a)^{1/2}g_2(bb^*)^{1/2})\,.
\end{equation*}
For details on operator spaces and completely bounded maps we refer
to the monographs \cite{ER} and \cite{Pi1}.

\section{Proof of the Effros-Ruan Conjecture} \setcounter{equation}{0}

%We first recall some preliminaries on Powers factors which will be
%needed for the proof of our main result.
Let $0< \lambda< 1$ be
fixed, and let $(\M\,, \phi)$ be the Powers factor of type
III$_\lambda$ with product state $\phi$\,, that is,
\[ (\M\,, \phi)=\bigotimes_{n=1}^\infty (M_2(\mathbb{C})\,,
\omega_\lambda)\,, \] where $\phi=\bigotimes _{n=1}^\infty
\omega_\lambda$\,, $\omega_\lambda(\,\cdot\,)=\text{Tr}(h_\lambda\,
\cdot\,)$ and $h_\lambda=\left(
\begin{array}
[c]{cc}%
\frac{\lambda}{1+\lambda} & 0\\
0 & \frac1{1+\lambda}
\end{array}
\right)$ (cf.\ \cite{Con})\,. The modular automorphism group
$(\sigma_t^{\phi})_{t\in \mathbb{R}}$ of $\phi$ is given by
\[ \sigma_t^{\phi}=\bigotimes_{n=1}^\infty
\sigma_t^{\omega_\lambda}\,, \] where for any matrix
$x=[x_{ij}]_{1\leq i, j\leq 2}\in M_2(\mathbb{C})$\,,
\[ \sigma_t^{\omega_\lambda}(x)=h_\lambda^{it}x h_\lambda^{-it}=\left(
\begin{array}
[c]{cc}%
x_{11} & \lambda^{it} x_{12}\\
\lambda^{-it} x_{21} & x_{22}
\end{array}
\right)\,, \quad t\in \mathbb{R}\,. \] Therefore
$\sigma_t^{\omega_\lambda}$ and $\sigma_t^{\phi}$ are periodic in
$t\in \mathbb{R}$ with minimal period
\[ t_0:=-\frac{2\pi}{\log \lambda}\,. \]
Let $\M_\phi$ denote the centralizer of $\phi$\,, that is,
\begin{equation*}
\M_\phi:=\{x\in \M: \sigma_t^\phi (x)=x\,, \forall t\in
\mathbb{R}\}\,.
\end{equation*}
It was proved by Connes (cf.\ \cite{Co4}\,, Theorem 4.26) that the
relative commutant of $\M_\phi$ in $\M$ is trivial, i.e.,
\[ \M_\phi^{\prime} \cap \M=\mathbb{C}1_{\M}\,, \]
where $1_{\M}$ denotes the identity of $\M$\,. In particular, $\phi$
is homogeneous in the sense of Takesaki (cf.\ \cite{Ta1}).

Furthermore, it is shown in \cite{Ha1}\,, (see Theorem 3.1 therein)
that the following strong Dixmier property holds for the Powers
factor $\M$\,. Namely, for all $x\in \M$\,,
\[ \phi(x)\cdot 1_{\M}\in \overline{\text{conv}\{vxv^*: v\in {\mathcal
U}(\M_\phi)\}}^{\|\cdot\|}\,, \] where the closure is taken in norm
topology and ${\mathcal U}(\M_\phi)$ denotes the unitary group on
$\M_\phi$\,. Moreover, by Corollary 3.4 in \cite{Ha1}\,, this can be
extended to finite sets in $\M$\,, i.e., for every finite set
$\{x_1\,, \ldots x_n\}\in \M$ and every $\varepsilon > 0$\,, there
exists a convex combination $\alpha$ of elements from
$\{\text{ad}(v): v\in {\mathcal U}(\M_\phi)\}$ such that
\begin{equation*}
\|\alpha(x_i)-\phi(x_i)\cdot 1_{\M}\|< \varepsilon\,, \quad 1\leq
i\leq n\,.
\end{equation*}
By standard arguments, it follows that there exists a net
$\{\alpha_i\}_{i\in I}\subseteq \text{conv}\{\text{ad}(v): v\in
{\mathcal U}(\M_\phi)\}$ such that
\begin{equation}\label{eq898908}
\lim\limits_{i\in I}\|\alpha_i(x)-\phi(x)\cdot 1_{\M}\|=0\,, \quad
x\in \M\,.
\end{equation}

In the following, we will identify $\M$ with $\pi_\phi(\M)$\,, where
$(\pi_\phi\,, H_\phi\,, \xi_\phi)$ is the GNS representation of $\M$
associated to the state $\phi$\,. Then
\[ H_\phi:=\overline{\M\xi_\phi}=L^2(\M\,, \phi)\,. \]
By Tomita-Takesaki theory (cf.\ \cite{Tk}), the operator $S_0$
defined by
\[ S_0(x\xi_\phi)=x^*\xi_\phi\,, \quad x\in \M \]
is closable. Its closure $S:=\overline{S_0}$ has a unique polar
decomposition
\begin{equation}\label{eq66666757555} S=J\Delta^{1/2}\,,
\end{equation} where $\Delta$ is a positive self-adjoint unbounded operator on
$L^2(\M\,, \phi)$ and $J$ is a conjugate-linear involution.
Moreover, for all $t\in \mathbb{R}$\,,
\[ \sigma_t^\phi (x)=\Delta^{it}x\Delta^{-it}\,, \quad x\in \M \] and
\[ J\M J=\M^\prime\,, \]
where $\M^\prime$ denotes the commutant of $\M$\,.

Following Takesaki's construction from \cite{Ta1}\,, define for all
$n\in \mathbb{Z}$
\[ \M_n:=\{x\in \M:\sigma_t^{\phi}(x)=\lambda^{int}x\,,\,\,\forall t\in
\mathbb{R}\}\,. \] Then, by Lemma 1.16 in \cite{Ta1}\,,
\[ \M_n=\{x\in
\M:\phi(xy)=\lambda^n\phi(yx)\,, \,\,\forall y\in \M\}\,. \] In
particular, $\M_{\phi}=\M_0$\,. It was proved in \cite{Ta1} (cf.\
Lemma 1.10) that $\M_n\neq \{0\}$\,, for all $n\in \mathbb{Z}$\,.
Furthermore, by a combination of Lemma 1.4 and Corollary 1.16 in
\cite{Ta1}\,, it follows that for all $n\in \mathbb{Z}$\,,
\begin{equation*}\Delta(\eta)=\lambda^n \eta\,, \quad \eta\in
\overline{\M_n\xi_\phi}
\end{equation*} and that
\[ L^2(\M\,, \phi)=\bigoplus\limits_{n=-\infty}^\infty
\overline{\M_n\xi_\phi}\,. \] As a consequence, one has the
following

\begin{lemma}\label{lem6786}
For every $n\in \mathbb{Z}$\,, there exists $c_n\in \M$ such that
\begin{equation}\label{eq67898}
\phi(c^*_nc_n)=\lambda^{-{n/2}}\,, \quad
\phi(c_nc_n^*)=\lambda^{n/2}
\end{equation}
and, moreover,
\begin{equation}\label{eq67899}
\langle c_nJc_nJ\xi_{\phi}\,, \xi_{\phi}\rangle_{H_{\phi}}=1\,.
\end{equation}
\end{lemma}

\begin{proof}
Let $n\in \mathbb{Z}$\,. Take $z\in \M_n\setminus \{0\}$\,. Then
$\phi(zz^*)=\lambda^n \phi(z^*z)$\,. Moreover,
\[
JzJ\xi_{\phi}=S\Delta^{-{1/2}}z\xi_{\phi}=S(\lambda^{-{n/2}}z)\xi_{\phi}=\lambda^{-{n/2}}z^*\xi_{\phi}\,.
\]
Therefore, $\langle zJzJ\xi_{\phi}\,,
\xi_{\phi}\rangle_{H_{\phi}}=\lambda^{-{n/2}}\langle
zz^*\xi_{\phi}\,,
\xi_{\phi}\rangle_{H_{\phi}}=\lambda^{-{n/2}}\phi(zz^*)=\lambda^{n/2}\phi(z^*z)$\,.
Hence \[ c_n:=(\lambda^{n/2}\phi(z^*z))^{-{1/2}}z \] satisfies
relations (\ref{eq67898}) and (\ref{eq67899})\,.
\end{proof}

Since $\M$ is an injective factor, it is known (cf.\ \cite{Con})
that for all finite sequences $x_1\,, \ldots \,,x_n\in \M$ and
$y_1\,, \ldots \,,y_n\in \M^\prime$\,, where $n$ is a positive
integer, the following holds
\begin{equation}\label{eq678992}
\left\|\sum\limits_{i=1}^n c_i d_i\right\|_{{\mathcal B}(L^2(\M\,,
\phi))}=\left\|\sum\limits_{i=1}^n c_i\otimes
d_i\right\|_{\M\otimes_{\text{min}} \M^\prime}\,.
\end{equation}
That is, the map defined by $c\otimes d\mapsto cd$\,, where $c\in
\M$ and $d\in \M^\prime$ extends uniquely to a C$^*$-algebra
isomorphism of
$\M\otimes_{\text{min}} \M^\prime$ onto $C^*(\M\,, \M^\prime)$\,.\\

Now let $A$ and $B$ be C$^*$-algebras and let $u:A\times
B\rightarrow \mathbb{C}$ be a jointly completely bounded bilinear
form.

\begin{prop}\label{prop7439}
There exists a bounded bilinear form
$\widehat{u}:(A\otimes_{\text{min}} \M)\times (B\otimes_{\text{min}}
\M^\prime)\rightarrow \mathbb{C}$ such that
\begin{equation}\label{eq4443444}\widehat{u}(a\otimes c\,, b\otimes d)=u(a\,, b)\langle
cd\xi_{\phi}\,, \xi_\phi \rangle_{H_\phi}\,, \quad a\in A\,, b\in
B\,, c\in \M\,, d\in {\M^\prime}\,,
\end{equation} and, moreover,
\begin{equation}\label{eq11111333331213}
\|\widehat{u}\|\leq \|u\|_{\text{jcb}}\,.
\end{equation}
\end{prop}

\begin{proof}
Let $a_1\,, \ldots\,, a_m\in A$\,, $b_1\,, \ldots \,, b_n\in B$\,,
$c_1\,, \ldots \,, c_m\in \M$\,, $d_1\,, \ldots \,, d_n\in
{\M^\prime}$\,, where $m$ and $n$ are positive integers. Then, by
(\ref{eq678992})\,,
\begin{eqnarray*}
\left|\sum_{i=1}^m \sum_{j=1}^n  u(a_i\,, b_j)\langle c_id_j
\xi_{\phi}\,, \xi_{\phi} \rangle_{H_{\phi}} \right| &\leq &
\left\|\sum_{i=1}^m \sum_{j=1}^n u(a_i\,, b_j)
c_i d_j\right\|_{{\mathcal B}(L^2(\M\,, \phi))}\\
&=& \left\|\sum_{i=1}^m \sum_{j=1}^n u(a_i\,,
b_j)c_i\otimes d_j\right\|_{\M\otimes_{\text{min}} \M^\prime}\\
&\leq & \|u\|_{\text{jcb}}\left\|\sum_{i=1}^m a_i\otimes
c_i\right\|_{A\otimes_{\text{min}} \M}\left\|\sum_{j=1}^n b_j\otimes
d_j\right\|_{B\otimes_{\text{min}} {\M^\prime}}\,,
\end{eqnarray*}
which yields the conclusion.
\end{proof}

\begin{lemma}\label{lem787878767678}
Let $v\in {\mathcal U}(\M_{\phi})$ and set $v^\prime=JvJ\in
{\M}^\prime$\,. Then, for all $x\in A\otimes \M$ and $y\in B\otimes
{\M}^\prime$\,,
\begin{equation}\label{eq454545434343}
\widehat{u}((\text{Id}_A\otimes \text{ad}(v))(x)\,,
(\text{Id}_B\otimes \text{ad}(v^\prime))(y))=\widehat{u}(x\,, y)\,.
\end{equation}
\end{lemma}

\begin{proof}
It suffices to prove that formula (\ref{eq454545434343}) holds for
elementary tensors $x=a\otimes c$ and $y=b\otimes d$\,, where $a\in
A$\,, $b\in B$\,, $c\in \M$ and $d\in {\M}^\prime$\,. By
(\ref{eq4443444}), it is enough to show that for all $c\in \M$ and
$d\in {\M}^\prime$\,,
\begin{equation}\label{eq33333222222323}
\langle vcv^*v^\prime d (v^\prime)^*)\xi_{\phi}\,,
\xi_{\phi}\rangle_{H_{\phi}}= \langle cd \xi_{\phi}\,,
\xi_{\phi}\rangle_{H_{\phi}}\,.
\end{equation}
Since $\{v\,, c\,, v^*\}$ commutes with $\{v^\prime\,, d\,,
(v^\prime)^*\}$\,, we have
\begin{eqnarray}\label{eq121211112121}
\langle vcv^*v^\prime d (v^\prime)^*\xi_{\phi}\,,
\xi_{\phi}\rangle_{H_{\phi}}&=& \langle v^\prime v cd v^*
(v^\prime)^*\xi_{\phi}\,, \xi_{\phi}\rangle_{H_{\phi}}\\
&=& \langle cdv^* (v^\prime)^*\xi_{\phi}\,,
v^*(v^\prime)^*\xi_{\phi}\rangle_{H_\phi}\,.\nonumber
\end{eqnarray}
But since $J\xi_{\phi}=\xi_{\phi}$\,, we deduce that
\begin{equation*}
v^*(v^\prime)^*\xi_{\phi}=v^*(JvJ)^*\xi_{\phi}=
v^*Jv^*J\xi_{\phi}=v^*Jv^*\xi_{\phi}\,.
\end{equation*}
Furthermore, since $v^*\in {\M}_{\phi}$ and $\Delta^{it}
\xi_{\phi}=\xi_{\phi}$\,, for all $t\in \mathbb{R}$\,, we have
\begin{equation*}
\Delta^{it}(v^*\xi_{\phi})=\sigma_t^{\phi}(v^*)\Delta^{it}\xi_{\phi}=v^*\xi_{\phi}\,,
\quad t\in \mathbb{R}\,.
\end{equation*}
Hence $v^*\xi_{\phi}$ is an eigenvector for $\Delta$ with
corresponding eigenvalue equal to 1. Using the polar decomposition
(\ref{eq66666757555}) of $S$\,, we infer that
\begin{equation*}
v^*Jv^*\xi_{\phi}=v^*S\Delta^{{-1}/2} u^*\xi_{\phi}=
v^*Sv^*\xi_{\phi}=v^*v\xi_{\phi}=\xi_{\phi}\,,
\end{equation*}
i.e., $v^*(v^{\prime})^*\xi_{\phi}=\xi_{\phi}$\,. Therefore
\begin{equation*}
\langle cd v^*(v^{\prime})^*\xi_{\phi}\,,
v^*(v^{\prime})^*\xi_{\phi}\rangle_{H_{\phi}}=\langle
cd\xi_{\phi}\,, \xi_{\phi}\rangle_{H_{\phi}}\,.
\end{equation*}
This gives (\ref{eq121211112121})\,, which completes the proof of
the lemma.
\end{proof}

\begin{lemma}\label{lem7876543}
Let $\{\alpha_i\}_{i\in I}\subseteq \text{conv}\{\text{ad}(v): v\in
{\mathcal U}(\M_\phi)\}$ be a net satisfying (\ref{eq898908})\,. For
every $i\in I$\,, consider the corresponding map $\alpha_i^\prime$
on $\M^\prime=J\M J$ given by
\begin{equation*}
\alpha_i^\prime (JxJ)=J\alpha_i(x)J\,, \quad x\in \M\,.
\end{equation*}
Moreover, let $\phi^\prime$ be the state on $\M^\prime$ defined by
\begin{equation}\label{eq8888898880}
\phi^\prime (JxJ):=\overline{\phi(x)}\,, \quad x\in \M\,.
\end{equation}
Furthermore, let $\hat{f}$ be a state on $A\otimes_{\text{min}} \M$
and $\hat{g}$ be a state on $B\otimes_{\text{min}} {\M^\prime}$\,,
arbitrarily chosen, and define states $f$ on $A$\,, respectively,
$g$ on $B$ by
\begin{eqnarray}
f(a)&=& \hat{f}(a\otimes 1_{\M})\,, \quad a\in A\label{eq7775789}\\
g(b)&=& \hat{g}(b\otimes 1_{\M^\prime})\,, \quad b\in
B\,,\label{eq77757893}
\end{eqnarray}
where $1_{\M^\prime}$ denotes the identity of $\M^\prime$\,. Then,
\begin{equation}\label{eq33334333}
\lim\limits_{i\in I} \hat{f}((\text{Id}_A \otimes
\alpha_i)(z))=(f\otimes \phi)(z)\,, \quad z\in A\otimes_{\text{min}}
\M\,,
\end{equation}
and, respectively,
\begin{equation}\label{eq333343334}
\lim\limits_{i\in I} \hat{g}((\text{Id}_B \otimes
\alpha_i^\prime)(w))=(g\otimes \phi^\prime)(w)\,, \quad w\in
B\otimes_{\text{min}} {\M}^\prime\,.
\end{equation}
\end{lemma}

\begin{proof}
Note that for $i\in I$\,, $\|\alpha_i\|_{\text{cb}}\leq 1$ and
$\|\alpha_i^\prime\|_{\text{cb}}\leq 1$\,. Therefore,
$\text{Id}_A\otimes \alpha_i$ and $\text{Id}_B\otimes
\alpha_i^\prime$ are well-defined contractions on
$A\otimes_{\text{min}} \M$ and $B\otimes_{\text{min}}
{\M}^\prime$\,, respectively. Hence, in order to prove
(\ref{eq33334333}) and (\ref{eq333343334})\,, it suffices to
consider elementary tensors $z=a\otimes c$ and $w=b\otimes d$\,,
where $a\in A$\,, $b\in B$\,, $c\in \M$ and $d\in \M^\prime$\,.

Let $a\in A$ and $c\in \M$\,. By (\ref{eq898908}) we deduce that the
following holds in norm topology \begin{equation*}\lim_{i\in
I}(\text{Id}_A \otimes \alpha_i)(a\otimes c)=\lim_{i\in I} a\otimes
\alpha_i(c)=\phi(c)(a\otimes 1_{\M})\,.\end{equation*} It follows
that
\[ \lim_{i\in I} \hat{f}((\text{Id}_A \otimes \alpha_i)(a\otimes
c))=\phi(c) \hat{f}(a\otimes 1_{\M})=\phi(c)f(a)=(f\otimes
\phi)(a\otimes c)\,,
\] which proves (\ref{eq33334333})\,. Further,
for all $x\in \M$\,,
\begin{eqnarray*}
\lim\limits_{i\in I} \alpha_i^\prime (JxJ)\,\,=\,\,\lim\limits_{i\in
I} J
\alpha_i(x) J &=& J(\phi(x)\cdot 1_{\M})J\\
&=&\overline{\phi(x)}J\cdot J \,\,=\,\,\overline{\phi(x)}\cdot
1_{\M}\,\,=\,\,\phi^\prime (JxJ)\cdot 1_{\M}\,,
\end{eqnarray*}
where the limit is taken in norm topology. Then (\ref{eq333343334})
can be proved in the same way as (\ref{eq33334333})\,.
\end{proof}

\begin{prop}\label{prop787654378}
Let $u$\,, $\widehat{u}$ and $\phi^\prime$ be as above. Then there
exist states $f_1\,, f_2$ on $A$ and states $g_1\,, g_2$ on $B$ such
that for all $x\in {A\otimes_{\text{min}} \M}$ and $y\in
{B\otimes_{\text{min}} {\M^\prime}}$\,,
\begin{equation}\label{eq6666678698}
|\widehat{u}(x, y)|\leq  \|u\|_{\text{jcb}}\left((f_1\otimes
\phi)(xx^*)+(f_2\otimes \phi)(x^*x)\right)^{1/2}\left((g_1\otimes
\phi^\prime)(y^*y)+(g_2\otimes \phi^\prime)(yy^*)\right)^{1/2}\,.
\end{equation}
\end{prop}

\begin{proof}
By the Grothendieck inequality for C$^*$-algebras (cf.\ \cite{Ha6})
applied to the bilinear form $\widehat{u}$\,, there exist states
$\hat{f_1}$\,, $\hat{f_2}$ on $A\otimes_{\text{min}} \M$ and states
$\hat{g_1}$\,, $\hat{g_2}$ on $B\otimes_{\text{min}} {\M^\prime}$
such that for all $x\in {A\otimes_{\text{min}} \M}$ and $y\in
{B\otimes_{\text{min}} {\M^\prime}}$\,,
\begin{eqnarray}\label{eq66666786}
|\widehat{u}(x, y)|&\leq
&\|\widehat{u}\|(\hat{f_1}(xx^*)+\hat{f_2}(x^*x))^{1/2}(\hat{g_1}(y^*y)+\hat{g_2}(yy^*))^{1/2}\\
&\leq
&\|u\|_{\text{jcb}}(\hat{f_1}(xx^*)+\hat{f_2}(x^*x))^{1/2}(\hat{g_1}(y^*y)+\hat{g_2}(yy^*))^{1/2}\nonumber
\,,
\end{eqnarray}
wherein we have used inequality (\ref{eq11111333331213})\,.

Since $\sqrt{\alpha\beta}\leq {(\alpha+\beta)}/{2}$ for all
$\alpha\,, \beta \geq 0$\,, it follows that
\begin{eqnarray*}
|\widehat{u}(x, y)|&\leq &
{\frac12}\|u\|_{\text{jcb}}\left(\hat{f_1}(xx^*)+\hat{f_2}(x^*x)
+\hat{g_1}(y^*y)+\hat{g_2}(yy^*)\right)\,.
\end{eqnarray*}
For $i=1\,, 2$\,, let $f_i$ be the state on $A$ constructed from
$\hat{f_i}$ by formula (\ref{eq7775789}), and, respectively, let
$g_i$ be the state on $B$ constructed from $\hat{g_i}$ by formula
(\ref{eq77757893}). We show in the following that these are the
states we are looking for.

By Lemma \ref{lem787878767678}\,, we deduce that for all $v\in
{\mathcal U}(\M_{\phi})$ (and $v^\prime:=JvJ$\,, as defined
therein)\,,
\begin{eqnarray}\label{eq444444434232}
|\widehat{u}(x, y)|&\leq &
{\frac12}\|u\|_{\text{jcb}}\Big[\hat{f_1}((\text{Id}_A\otimes
\text{ad}(v))(xx^*))+\hat{f_2}((\text{Id}_A\otimes
\text{ad}(v))(x^*x))+ \\ &&\qquad \qquad  +
\hat{g_1}((\text{Id}_B\otimes
\text{ad}(v^\prime))(y^*y))+\hat{g_2}((\text{Id}_B\otimes
\text{ad}(v^\prime))(yy^*))\Big]\,.\nonumber
\end{eqnarray}
Next choose nets $\{\alpha_i\}_{i\in I}$ and
$\{\alpha_i^\prime\}_{i\in I}$ as in Lemma \ref{lem7876543}\,. For
all $i\in I$\,, it follows that
\begin{eqnarray}\label{eq4444444342327}
|\widehat{u}(x, y)|&\leq &
{\frac12}\|u\|_{\text{jcb}}\Big[\hat{f_1}((\text{Id}_A\otimes
\alpha_i)(xx^*))+\hat{f_2}((\text{Id}_A\otimes \alpha_i)(x^*x))+
\\ &&\qquad \qquad  + \hat{g_1}((\text{Id}_B\otimes
\alpha_i^\prime)(y^*y))+\hat{g_2}((\text{Id}_B\otimes
\alpha_i^\prime)(yy^*))\Big]\,,\nonumber
\end{eqnarray}
since the right-hand side of (\ref{eq4444444342327}) is a convex
combination of the possible right-hand sides of
(\ref{eq444444434232})\,. Then, by Lemma \ref{lem7876543} we obtain
in the limit that
\begin{equation}\label{eq666667867332}
|\widehat{u}(x, y)|\leq  {\frac12}\|u\|_{\text{jcb}}((f_1\otimes
\phi)(xx^*)+(f_2\otimes \phi)(x^*x) + (g_1\otimes
\phi^\prime)(y^*y)+(g_2\otimes \phi^\prime)(yy^*))\,.
\end{equation}
Recall that $x$ and $y$ were arbitrarily chosen in
${A\otimes_{\text{min}} \M}$ and ${B\otimes_{\text{min}}
{\M^\prime}}$\,, respectively. Hence, replacing $x$ by $t^{1/2}x$
and $y$ by $t^{{-1}/{2}}y$\,, where $t> 0$\,, we deduce that the
following inequality holds for all $x\in {A\otimes_{\text{min}}
\M}$\,, $y\in {B\otimes_{\text{min}} {\M^\prime}}$ and $t> 0$\,:
\begin{equation}\label{eq66666786733255}
|\widehat{u}(x, y)|\leq {\frac12}
\|u\|_{\text{jcb}}\left(t(f_1\otimes \phi)(xx^*)+t(f_2\otimes
\phi)(x^*x) + {\frac1t}(g_1\otimes
\phi^\prime)(y^*y)+{\frac1t}(g_2\otimes \phi^\prime)(yy^*)\right)\,.
\end{equation}
Since for all $\alpha\,, \beta \geq 0$\,, we have
\begin{equation}\label{eq22222452}
\inf_{t> 0} (t\alpha+t^{-1}\beta)=2\sqrt{\alpha\beta}\,,
\end{equation}
the assertion then follows by taking infimum over all $t> 0$ in
(\ref{eq66666786733255}).
\end{proof}

\begin{lemma}\label{lem44443447}
Let $\alpha\,, \beta\geq 0$\,. Then
\begin{equation}\label{eq343434546}
\inf_{n\in \mathbb{Z}} (\lambda^n \alpha+\lambda^{-n} \beta)\leq
(\lambda^{1/2}+\lambda^{{-1}/2})\sqrt{\alpha\beta}\,.
\end{equation}
\end{lemma}

\begin{proof}
The statement is obvious if $\alpha=0$ or $\beta=0$\,. Assume that
$\alpha\,, \beta> 0$\,. Since $0< \lambda< 1$\,, it follows that
$(0\,, \infty)=\bigcup_{n\in \mathbb{Z}} [\lambda^{2n+1}\,,
\lambda^{2n-1}]$\,. Hence, we can choose $n\in \mathbb{Z}$ such that
\begin{equation*}\label{eq121211111}
\lambda^{2n+1}\leq {\beta}/{\alpha}\leq \lambda^{2n-1}\,.
\end{equation*}
Set $\alpha_1:=\lambda^n \alpha$ and $\beta_1:=\lambda^{-n}
\beta$\,. Then $\lambda\leq {\beta_1}/{\alpha_1}\leq {1/\lambda}$\,.
Since the function $t\mapsto t^{1/2}+t^{{-1}/2}$ is decreasing on
$[\lambda\,, 1]$ and increasing on $[1\,, {1/{\lambda}}]$\,, it
follows that
\begin{equation*}
\max \{t^{1/2}+t^{{-1}/2}: t\in [\lambda\,, 1/{\lambda}]\}=
\lambda^{1/2}+\lambda^{{-1}/2}\,.
\end{equation*}
Hence, we deduce that
\begin{eqnarray*}
\lambda^n \alpha+ \lambda^{-n} \beta\,\,=\,\,\alpha_1+\beta_1 &=&
\left(\sqrt{{\alpha_1}/{\beta_1}}+\sqrt{{\beta_1}/{\alpha_1}}\right)\sqrt{\alpha_1
\beta_1}\\&\leq &(\lambda^{1/2}+\lambda^{{-1}/2})\sqrt{\alpha_1
\beta_1}\,\,=\,\,(\lambda^{1/2}+\lambda^{{-1}/2})\sqrt{\alpha
\beta}\,,
\end{eqnarray*}
which proves (\ref{eq343434546})\,.
\end{proof}

\begin{prop}\label{prop89999998} Set
\[C(\lambda):=\sqrt{{\left(\lambda^{1/2}+\lambda^{-1/2}\right)}/{2}}\,.
\] Let $u$ be as above and let $f_1\,, f_2$ be states on $A$\,, respectively, $g_1\,, g_2$
be states on $B$ as in Proposition \ref{prop787654378}\,. Then, for
all $a\in A$ and $b\in B$\,,
\begin{equation}\label{eq555555786}
|u(a, b)|\leq
C(\lambda)\|u\|_{\text{jcb}}\left(f_1(aa^*)^{1/2}g_1(b^*b)^{1/2}
+f_2(a^*a)^{1/2}g_2(bb^*)^{1/2}\right)\,.
\end{equation}
that is, the Effros-Ruan conjecture holds with constant
$C(\lambda)$\,.
\end{prop}

\begin{proof}
Let $n\in \mathbb{Z}$ and choose $c_n\in \M$ as in Lemma
\ref{lem6786}\,. Then, for all $a\in A$ and all $b\in B$\,, it
follows by (\ref{eq4443444}) and (\ref{eq67899}) that
\begin{equation*}
\widehat{u}(a\otimes c_n\,, b\otimes J{c_n}J)= u(a, b)\langle
{c_n}J{c_n}J\xi_{\phi}\,, \xi_{\phi}\rangle_{H_{\phi}}= u(a, b)\,.
\end{equation*}
By Proposition \ref{prop787654378}, together with
(\ref{eq8888898880}) and (\ref{eq67898})\,, it follows that
\begin{eqnarray}\label{eq7777776777890}
\hspace{1cm}|u(a, b)|^2&=& |\widehat{u}(a\otimes {c_n}\,, b\otimes
{J{c_n}J})|^2\\
&\leq
&\|u\|^2_{\text{jcb}}\left(f_1(aa^*)\phi({c_n}c_n^*)+f_2(a^*a)\phi(c_n^*{c_n})\right)\left(g_1(b^*b)\phi({c_n}^*{c_n})+g_2(bb^*)\phi({c_n}{c_n}^*)\right)\nonumber\\
&=&
\|u\|^2_{\text{jcb}}\left(\lambda^{{n}/2}f_1(aa^*)+\lambda^{{-n}/2}
f_2(a^*a)\right)\left(\lambda^{{-n}/2}
g_1(b^*b)+\lambda^{{n}/2}g_2(bb^*)\right)\nonumber\\
&=&
\|u\|^2_{\text{jcb}}\left(f_1(aa^*)g_1(b^*b)+f_2(a^*a)g_2(bb^*)+\lambda^{n}
f_1(aa^*)g_2(bb^*)+\lambda^{-n}
f_2(a^*a)g_1(b^*b)\right)\,.\nonumber
\end{eqnarray}
Note that $\lambda^{1/2}+\lambda^{{-1}/2}=2C(\lambda)^2$\,. By
taking infimum in (\ref{eq7777776777890}) over all $n\in
\mathbb{Z}$\,, we deduce from Lemma \ref{lem44443447} that
\begin{eqnarray*}
|u(a, b)|^2&\leq &
\|u\|^2_{\text{jcb}}\left(f_1(aa^*)g_1(b^*b)+f_2(a^*a)g_2(bb^*)+2C(\lambda)^2f_1(a^*a)^{\frac12}g_1(b^*b)^{\frac12}f_2(aa^*)^{\frac12}g_2(bb^*)^{\frac12}\right)\\
&\leq &
C(\lambda)^2\|u\|^2_{\text{jcb}}\left(f_1(aa^*)^{\frac12}g_1(b^*b)^{\frac12}+f_2(a^*a)^{\frac12}g_2(bb^*)^{\frac12}\right)^2\,,
\end{eqnarray*}
wherein we have used the fact that $C(\lambda)> 1$\,. The assertion
follows now by taking square roots.
\end{proof}

\noindent{\em Proof of Theorem \ref{theffruan}}:
 Thus far we have proved that given C$^*$-algebras
$A$ and $B$ and a j.c.b. bilinear form $u:A\times B\rightarrow
\mathbb{C}$\,, then the Effros-Ruan conjecture holds with constant
$C(\lambda)=\sqrt{{\left(\lambda^{1/2}+\lambda^{-1/2}\right)}/{2}}$\,,
for every $0< \lambda< 1$\,. Now recall that the sets
\begin{equation*}
Q(A):=\{f\in A_{+}^*: \|f\|\leq 1\}\,, \quad Q(B):=\{g\in B_{+}^*:
\|g\|\leq 1\}
\end{equation*}
are compact in the weak$^*$-topology\,, where $A_{+}^*$ and
$B_{+}^*$ denote the sets of positive functionals on $A$ and $B$\,,
respectively. Since $C(\lambda)\rightarrow 1$ as $\lambda\rightarrow
1$\,, by using a simple compactness argument it follows from
Proposition \ref{prop89999998} that there exist $f_1^0\,, f_2^0\in
Q(A)$ and $g_1^0\,, g_2^0\in Q(B)$ such that for all $a\in A$ and
$b\in B$\,,
\begin{equation*}
|u(a, b)|\leq
\|u\|_{\text{jcb}}\left(f_1^0(aa^*)^{1/2}g_1^0(b^*b)^{1/2}
+f_2^0(a^*a)^{1/2}g_2^0(bb^*)^{1/2}\right)\,.
\end{equation*}
But $f_i^0\leq f_i$\,, respectively, $g_i^0\leq g_i$\,, $i=1, 2$\,,
where $f_1\,, f_2$ are states on $A$ and $g_1\,, g_2$ are states on
$B$\,. Therefore the Effros-Ruan conjecture holds with constant
one.\qed

\section{Applications}
\setcounter{equation}{0}

Let $E\subseteq A$ and $F\subseteq B$ be operator spaces sitting in
C$^*$-algebras $A$ and $B$\,. Let $u:E\times F\rightarrow
\mathbb{C}$ be a bounded bilinear form. Define $\|u\|_{\text{ER}}$
to be the smallest constant $0\leq \kappa_2\leq \infty$ for which
there exist states $f_1$\,, $f_2$ on $A$ and states $g_1$\,, $g_2$
on $B$ such that for all $a\in E$ and $b\in F$\,,
\begin{equation}\label{eq667766776767}
|u(a, b)|\leq \kappa_2
(f_1(aa^*)^{1/2}g_1(b^*b)^{1/2}+f_2(a^*a)^{1/2}g_2(bb^*)^{1/2})\,.
\end{equation}
In the case when $E=A$ and $F=B$\,, we have from Theorem
\ref{theffruan} that $\|u\|_{\text{ER}}\leq \|u\|_{\text{jcb}}$\,.
Moreover, if $E$ and $F$ are exact operator spaces and $u:E\times
F\rightarrow \mathbb{C}$ is a j.c.b.\ bilinear form, then by
\cite{PiSh} (cf.\ Theorem 0.3 and 0.4)\,,
\begin{equation*}
\|u\|_{\text{ER}}\leq
2^{3/2}\text{ex}(E)\text{ex}(F)\|u\|_{\text{jcb}}\,.
\end{equation*}
However, for bilinear forms on general operator spaces $E$ and $F$
it can happen that $\|u\|_{\text{jcb}}< \infty$\,, while
$\|u\|_{\text{ER}}=\infty$ (see Example \ref{ex878765456} below)\,.
Therefore Theorem \ref{theffruan} cannot be generalized to arbitrary
operator spaces.

Recall that a bilinear map $u:E\times F\rightarrow \mathbb{C}$ is
called {\em completely bounded} (in the sense of Christensen and
Sinclair) (see \cite{ChS}\,, \cite{PiSh} and the references given
therein) if the bilinear forms $u_n:M_n(E)\times M_n(F)\rightarrow
M_n(\mathbb{C})$ defined by
\begin{equation*}
u_n(a\otimes x\,, b\otimes y):=u(a\,, b)xy\,, \quad a\in E\,, b\in
F\,, x\,, y\in M_n(\mathbb{C})
\end{equation*}
are uniformly bounded, in which case we set
\begin{equation}\label{eq888565658}
\|u\|_{\text{cb}}:=\sup\limits_{n\in \mathbb{N}} \|u_n\|\,.
\end{equation}
Moreover, $u$ is completely bounded if and only if there exists a
constant $\kappa_3\geq 0$ and states $f$ on $A$ and $g$ on $B$ such
that for all $a\in E$ and $b\in F$\,,
\begin{equation}\label{eq33333335554555}
|u(a, b)|\leq \kappa_3f(aa^*)^{1/2}g(b^*b)^{1/2}
\end{equation}
and $\|u\|_{\text{cb}}$ is the smallest constant $\kappa_3$ for
which (\ref{eq33333335554555}) holds (see also the Introduction to
\cite{PiSh}).

It was shown by Effros and Ruan (cf.\ \cite{ER3}) that if $u:E\times
F\rightarrow \mathbb{C}$ is completely bounded, then the associated
map $\widetilde{u}:E\rightarrow F^*$ defined by (\ref{eq555555567})
admits a factorization of the form $\widetilde{u}=vw$ through a row
Hilbert space $H_r$\,, where $E \overset{v}{\longrightarrow} H_r
\overset{w}{\longrightarrow} F^*$ and
$\|v\|_{\text{cb}}\|w\|_{\text{cb}}=\|u\|_{\text{cb}}$\,. In
particular, it follows that
\begin{equation}\label{1111132}
\|u\|_{\text{jcb}}:=\|\widetilde{u}\|_{\text{cb}}\leq
\|u\|_{\text{cb}}\,.
\end{equation}

\begin{lemma}\label{lem5675321232} {\rm{(}}cf.\ \cite{PiSh} and
\cite{Xu}{\rm{)}} Let $u:E\times F\rightarrow \mathbb{C}$ be a
bounded bilinear form on operator spaces $E\subseteq A$ and
$F\subseteq B$ sitting in C$^*$-algebras $A$ and $B$. Let $f_1$\,,
$f_2$ be states on $A$ and $g_1$\,, $g_2$ be states on $B$ such that
for all $a\in E$ and $b\in F$\,,
\begin{equation*}
|u(a, b)|\leq
\|u\|_{\text{ER}}(f_1(aa^*)^{1/2}g_1(b^*b)^{1/2}+f_2(a^*a)^{1/2}g_2(bb^*)^{1/2})\,.
\end{equation*}
Then $u$ can be decomposed as $u=u_1+ u_2$\,, where $u_1$ and $u_2$
are bilinear forms satisfying the following inequalities, for all
$a\in A$ and $b\in B$\,:
\begin{eqnarray}
|u_1(a, b)|&\leq
&\|u\|_{\text{ER}}f_1(aa^*)^{1/2}g_1(b^*b)^{1/2}\label{eq456}\\
|u_2(a, b)|&\leq &
\|u\|_{\text{ER}}f_2(a^*a)^{1/2}g_2(bb^*)^{1/2}\,.\label{eq457}
\end{eqnarray}
In particular,
\begin{equation*}
\|u_1\|_{\text{cb}}\leq \|u\|_{\text{ER}}\,, \quad
\|u_2^{t}\|_{\text{cb}}\leq \|u\|_{\text{ER}}\,,
\end{equation*}
where $u_2^{t}(b, a):=u_2(a, b)$\,, for all $a\in E$ and $b\in F$\,.
\end{lemma}

\begin{proof}
Such a decomposition was obtained in \cite{PiSh} (cf.\ last
statement in Theorem 0.4 in \cite{PiSh}), except that the states
$f_1$\,, $f_2$\,, $g_1$\,, $g_2$ satisfying (\ref{eq456}) and
(\ref{eq457}) were possibly different from the original ones. Later,
following a suggestion of Pisier, Xu proved the above decomposition
without change of states. (See \cite{Xu}, Proposition 5.1 and the
Remark following the proof of this proposition.)
\end{proof}

\begin{rem}\label{rem56743233332}\rm
Note that our main result combined with the above splitting lemma
solves conjecture $(0.2')$ in \cite{PiSh} (with constant $K=2$)\,,
and hence it solves Blecher's conjecture (cf.\ \cite{Bl} and
Conjecture $(0.2)$ in \cite{PiSh}).
\end{rem}

\begin{prop}\label{prop8989898999}
$(i)$ \,\,Let $u:A\times B\rightarrow \mathbb{C}$ be a bounded
bilinear form on C$^*$-algebras $A$ and $B$\,. Then
\begin{equation}\label{eq472}
\|u\|_{ER}\leq \|u\|_{jcb}\leq 2\|u\|_{ER}\,.
\end{equation}
$(ii)$ \,\,Let $c_1\,, c_2$ denote the best constants in the
inequalities
\begin{equation}\label{eq4724724472}
c_1\|u\|_{ER}\leq \|u\|_{jcb}\leq c_2\|u\|_{ER}\,,
\end{equation}
where $u:A\times B\rightarrow \mathbb{C}$ is any bounded bilinear
form on arbitrary C$^*$-algebras $A$ and $B$\,. Then $c_1=1$ and
$c_2=2$\,.
\end{prop}

\begin{proof}
$(i)$\,. The left-hand side inequality follows from our main
theorem, while the right-hand side inequality follows from the
splitting lemma above. Indeed, we can assume that
$\|u\|_{\text{ER}}< \infty$\,. Then with $u_1\,, u_2:A\times
B\rightarrow \mathbb{C}$ as in Lemma \ref{lem5675321232}\,,
\begin{equation*}
\|u\|_{\text{jcb}}\leq \|u_1\|_{\text{jcb}}+\|u_2\|_{\text{jcb}}
=\|u_1\|_{\text{jcb}}+\|u_2^{t}\|_{\text{jcb}}\leq
\|u_1\|_{\text{cb}}+\|u_2^{t}\|_{\text{cb}}\leq
2\|u\|_{\text{ER}}\,.
\end{equation*}
$(ii)$\,. By $(i)$ we know that $c_1\geq 1$ and $c_2\leq 2$\,. We
now prove that $c_2=2$\,. Let $\tau$ be a tracial state on a
C$^*$-algebra $A$ and define a bilinear form $u:A\times A\rightarrow
\mathbb{C}$ by $u(a, b):=\tau(ab)$\,, for all $a\,, b\in A$\,. Then
$\|u\|_{\text{jcb}}\geq \|u\|=1$\,, and for all $a\,, b\in A$\,,
\begin{equation*}
|u(a, b)|\leq \tau(aa^*)^{1/2}\tau(b^*b)^{1/2}
=\frac12\left(\tau(aa^*)^{1/2}\tau(b^*b)^{1/2}+\tau(a^*a)^{1/2}\tau(bb^*)^{1/2}\right)\,,
\end{equation*}
which implies that $\|u\|_{ER}\leq {\frac12}$. By (\ref{eq472})\,,
$\|u\|_{ER}\geq {\frac12} \|u\|_{\text{jcb}}$\,. Hence
$\|u\|_{ER}=\frac12$ and $\|u\|_{\text{jcb}}=1$\,, and the assertion
follows. To prove that $c_1=1$\,, let $\phi$ be any state on a
unital, properly infinite C$^*$-algebra $A$\,. Let $u:A\otimes
A\rightarrow \mathbb{C}$  be defined by $u(a, b):=\phi(ab)$\,, for
all $a\,, b\in A$\,. Note that
\[ \|u\|_{ER}\leq \|u\|_{\text{jcb}}\leq \|u\|_{\text{cb}}\leq 1\,,
\]
where the last inequality follows immediately from (${0.5}^\prime$)
in \cite{PiSh} (by taking $f_1=g_1=\phi$ therein). We claim that
$\|u\|_{ER}=1$\,. For this, let $f_1\,, f_2\,, g_1\,, g_2$ be states
on $A$ and let $\{s_n\}_{n\geq 1}$ be a sequence of isometries in
$A$ with orthogonal ranges. Then $f_k(s_n s_n^*)\rightarrow 0$ as
$n\rightarrow \infty$\,, respectively $g_k(s_n s_n^*)\rightarrow 0$
as $n\rightarrow \infty$\,, for $k=1\,, 2$\,. Note that $u(s_n\,,
s_n^*)=1$\,, for all $n\geq 1$\,, while
\[ \lim_{n\rightarrow \infty} f_1(s_n s_n^*)^{1/2} g_1(s_n s_n^*)^{1/2}+ f_2(s_n^*
s_n)^{1/2}g_2(s_n^* s_n)^{1/2}=1\,. \] This shows that
$\|u\|_{ER}\geq 1$ and the assertion is proved.
\end{proof}

\begin{lemma}\label{lem6677667756765}
Let $E\subseteq A$ and $F\subseteq B$ be operator spaces sitting in
C$^*$-algebras $A$ and $B$\,, and let $u:E\times F\rightarrow
\mathbb{C}$ be a bounded bilinear form. If $\|u\|_{ER}<\infty$\,,
then the associated map $\widetilde{u}:E\rightarrow F^*$ admits a
cb-factorization $\widetilde{u}=vw$ through $H_r\oplus K_c$ for some
Hilbert spaces $H$ and $K$\,, where $E \overset{v}{\longrightarrow}
H_r\oplus K_c \overset{w}{\longrightarrow} F^*$\,, satisfying
\begin{equation*}
\|v\|_{cb}\|w\|_{cb}\leq 2\|\widetilde{u}\|_{ER}\,.
\end{equation*}
%$$
%\xymatrix{
% {E}\ar@{->}^{\widetilde{u}}[rr]
% \ar@{->}_{v}[dr]
% & & {F^*}\\
%& {H_r\oplus K_c} \ar@{->}[ur]_{w}}
%$$
%with $\|v\|_{cb}\|w\|_{cb}\leq 2\|\widetilde{u}\|_{ER}$\,.
\end{lemma}

\begin{proof}
Choose states $f_1\,, f_2$ on $A$ and states $g_1\,, g_2$ on $B$
such that (\ref{eq667766776767}) holds. Then, by Lemma
\ref{lem5675321232}\,, $u$ can be decomposed as $u=u_1+u_2$\,, where
$u_1$ and $u_2$ are bounded bilinear forms satisfying (\ref{eq456})
and (\ref{eq457}). The rest of the proof follows from the proof of
Corollary 0.7 on p. 206 in \cite{PiSh}.
\end{proof}

\begin{prop}\label{prop565656445}
Let $A$ and $B$ be C$^*$-algebras. Then every completely bounded
linear map $T:A\rightarrow B^*$ admits a cb-factorization $T=vw$
through $H_r\oplus K_c$ for some Hilbert spaces $H$ and $K$\,, such
that
\begin{equation*}\|u\|_{\text{cb}}\|w\|_{\text{cb}}\leq
2\|T\|_{\text{cb}}\,.\end{equation*}
\end{prop}

\begin{proof}
Let $T:A\rightarrow B^*$ be a completely bounded linear map. Then
$T$ is of the form $T=\widetilde{u}$\,, for a j.c.b. bilinear form
$u:A\times B\rightarrow \mathbb{C}$ with
$\|u\|_{\text{jcb}}=\|T\|_{\text{cb}}$\,. The assertion follows now
from Lemma \ref{lem6677667756765}, by using the fact that
$\|u\|_{ER}\leq \|u\|_{\text{jcb}}$\,.
\end{proof}

The following example is implicit in the proof of Corollary 3.2 in
\cite{PiSh}:
\begin{exam}\label{ex878765456}\rm
Let $E$ be an operator space which is not Banach space isomorphic to
a Hilbert space, and let $E\subseteq A$ and $E^*\subseteq B$ be
completely isometric embeddings of $E$ and $E^*$\,, respectively,
into C$^*$-algebras $A$ and $B$\,. Define $u:E\times E^*\rightarrow
\mathbb{C}$ by
\begin{equation*}
u(a, b):=b(a)\,, \quad a\in E\,, b\in E^*\,.
\end{equation*}
Then $\widetilde{u}:E\rightarrow E^{**}$ is the standard inclusion
of $E$ into its second dual. Therefore
$\|u\|_{\text{jcb}}=\|\widetilde{u}\|_{\text{cb}}=1$\,. We will show
that $\|u\|_{\text{ER}}=\infty$\,. If $\|u\|_{\text{ER}}< \infty$\,,
then it follows from Lemma \ref{lem6677667756765} that
$\widetilde{u}$ admits a cb-factorization through $H_r\oplus K_c$\,,
for some Hilbert spaces $H$ and $K$\,. In particular,
$\widetilde{u}:E\rightarrow E^{**}$ has a Banach space factorization
through a Hilbert space. This contradicts the assumption on $E$\,.
Hence $\|u\|_{\text{ER}}=\infty$\,.
\end{exam}

The following result was proved in \cite{PiSh} with constant
$2^{9/4}$ instead of $\sqrt{2}$ (see the second part of Corollary
3.4 in \cite{PiSh})\,.

\begin{cor}\label{cor7777787777}
Let $T$ be a completely bounded linear map from a C$^*$-algebra $A$
to the operator Hilbert space $OH(I)$\,, $I$ being an arbitrary
index set. Then there exist states $f_1$ and $f_2$ on $A$ such that
\begin{equation*}
\|T(a)\|\leq \sqrt{2} f_1(aa^*)^{1/4} f_2(a^*a)^{1/4}\,, \quad a\in
A\,.
\end{equation*}
\end{cor}

\begin{proof}
Given a vector space $E$\,, we let $\bar{E}$ denote the conjugate
vector space. Let $J:OH(I)\rightarrow \overline{{OH(I)}^*}$ be the
canonical cb-isomorphism of $OH(I)$ with the conjugate of its dual
space (cf. \cite{Pi6}), and set
\[ V:=\overline{T^*} J T\,, \]
where $T^*: {OH(I)}^*\rightarrow A^*$ is the adjoint of $T$\,. Then
$V$ is a completely bounded linear map from $A$ to
$\overline{A^*}=(\bar{A})^*$\,. Therefore $V=\widetilde{v}$ for a
j.c.b. bilinear form $v:A\times \bar{A}\rightarrow \mathbb{C}$\,.
Moreover,
\[ \|v\|_{\text{jcb}}=\|V\|_{\text{cb}}\leq \|T\|_{\text{cb}}^2\,.
\]
Actually, equality holds above (cf.\ \cite{PiSh}, proof of Corollary
3.4), but we shall not need this. By our main theorem, there exist
states $f_1^0\,, f_2^0$ on $A$ and states $g_1^0\,, g_2^0$ on
$\bar{A}$ such that for all $a\in A$ and $b\in \bar{A}$\,,
\begin{equation*}
|v(a, b)|\leq
\|T\|_{\text{cb}}^2\left(f_1^0(aa^*)^{1/2}g_1^0(b^*b)^{1/2}+f_2^0(a^*a)^{1/2}g_2^0(bb^*)^{1/2}\right)\,.
\end{equation*}
The canonical isomorphism $J$ of $OH(I)$ onto $\overline{{OH(I)}^*}$
satisfies
\[ \overline{J(x)}(x)=\|x\|^2=J(x)(\bar{x})\,, \quad x\in OH(I)\,.
\]
For all $a\in A$ we then have $v(a\,,
\bar{a})=(Va)(\bar{a})=(\overline{T^*} J Ta)(\bar{a})=(J T
a)(\,\overline{Ta}\,)=\|Ta\|^2$\,, and therefore
\begin{eqnarray*}
\|Ta\|^2=|v(a\,, \bar{a})|&\leq &
\|T\|_{\text{cb}}^2\left(f_1^0(aa^*)^{1/2}g_1^0(\,\overline{a^*a}\,)^{1/2}+f_2^0(a^*a)^{1/2}g_2^0(\,\overline{aa^*}\,)^{1/2}\right)\\
&\leq &
\|T\|_{\text{cb}}^2\left(f_1^0(aa^*)+g_2^0(\,\overline{aa^*}\,)\right)^{1/2}\left(f_2^0(a^*a)+g_1^0(\,\overline{a^*a}\,)\right)^{1/2}\\
&\leq & 2\|T\|_{\text{cb}}^2f_1(aa^*)^{1/2}f_2(a^*a)^{1/2}\,,
\end{eqnarray*}
where $f_1$ and $f_2$ are states on $A$ given by
\begin{equation*}
f_1(a):={\frac12}\left(f_1^0(a)+\overline{g_2^0(\bar{a})}\right)\,,
\quad
f_2(a):={\frac12}\left(f_2^0(a)+\overline{g_1^0(\bar{a})}\right)\,,
\quad a\in A\,.
\end{equation*}
This completes the proof.
\end{proof}

As a consequence of Proposition 3.3 we also obtain (by adjusting the
corresponding proofs in \cite{PiSh}) the following strengthening of
Corollaries 3.1 and 3.3 in \cite{PiSh}:

\begin{cor}\label{cor888899888}
Let $E$ be an operator space such that $E$ and its dual $E^*$ embed
completely isomorphically into preduals $M_*$ and $N_*$\,,
respectively, of von Neumann algebras $M$ and $N$\,. Then $E$ is
cb-isomorphic to a quotient of a subspace of $H_r\oplus K_c$\,, for
some Hilbert spaces $H$ and $K$\,.
\end{cor}

\begin{cor}\label{cor8888998885}
Let $E$ be an operator space and let $E\subseteq A$ and
$E^*\subseteq B$ be completely isometric embeddings into
C$^*$-algebras $A$ and $B$ such that both subspaces are completely
complemented. Then $E$ is cb-isomorphic to $H_r\oplus K_c$ for some
Hilbert spaces $H$ and $K$\,.
\end{cor}

Note that as another consequence of our main theorem we obtain (with
essentially the same proof as the corresponding Corollary 0.6 in
\cite{PiSh}) the following result:

\begin{cor}\label{cor66665655555}
Let $A_0$\,, $A$\,, $B_0$ and $B$ be C$^*$-algebras such that
$A_0\subseteq A$ and $B_0\subseteq B$\,. Then any j.c.b. bilinear
form $u_0:A_0\times B_0\rightarrow \mathbb{C}$ extends to a bilinear
form $u:A\times B\rightarrow \mathbb{C}$ such that
\begin{equation*}
\|u\|_{\text{jcb}}\leq 2\|u_0\|_{\text{jcb}}\,.
\end{equation*}
\end{cor}

Let $u:A\times B\rightarrow \mathbb{C}$ be a j.c.b. bilinear form on
C$^*$-algebras $A$ and $B$\,. Recall that $\|u\|_{\text{jcb}}$ is
the smallest constant $\kappa_1$ for which inequality
(\ref{eq3333222333220}) holds, for arbitrary C$^*$-algebras $C$ and
$D$\,. The following result shows that if the inequality
(\ref{eq3333222333220}) holds for the given bilinear form $u$ with
constant $\kappa_1$\,, then the same inequality (with $\kappa_1$
replaced by $2\kappa_1$) holds for $u$\,, when the
$(C\otimes_{\text{min}} D)$-norm on the left-hand side is replaced
by the $(C\otimes_{\text{max}} D)$-norm.

\vspace*{0.1cm}
\begin{prop}\label{prop789}
Let $A$ and $B$ be C$^*$-algebras, and let $u: A\times B\rightarrow
\mathbb{C}$ be a j.c.b.\ bilinear form. Then, for all C$^*$-algebras
$C$ and $D$\,, all $m, n\in \mathbb{N}$ and all finite sequences
$a_1\,, \ldots \,, a_m\in A$\,, $b_1\,, \ldots \,, b_n\in B$\,,
$c_1\,, \ldots \,, c_m\in C$\,, $d_1\,, \ldots \,, d_n\in D$\,,
\begin{equation}\label{eq222211112} \left\|\sum\limits_{i=1}^m \sum\limits_{j=1}^n u(a_i, b_j)c_i\otimes
d_j \right\|_{C\otimes_{\text{max}} D}\leq
2\|u\|_{\text{jcb}}\left\|\sum\limits_{i=1}^m a_i\otimes
c_i\right\|_{A\otimes_{\text{min}} C}\left\|\sum\limits_{j=1}^n
b_j\otimes d_j\right\|_{B\otimes_{\text{min}} D}\,.
\end{equation}
\end{prop}

\begin{proof}
There exist states $f_1\,, f_2$ on $A$ and $g_1\,, g_2$ on $B$ such
that inequality (\ref{eq667766776767}) holds. Then, as explained in
the proof of Lemma \ref{lem6677667756765}\,, $u$ can be decomposed
as $u=u_1+ u_2$\,, where $u_1$ and $u_2$ are bounded bilinear forms
satisfying (\ref{eq456}) and (\ref{eq457})\,.

By the definition of $\|\cdot\|_{\text{max}}$\,, in order to prove
(\ref{eq222211112}) we have to show that for all pairs of commuting
representations $\pi:A\rightarrow {\mathcal B}(H)$\,,
$\rho:B\rightarrow {\mathcal B}(H)$\,, where $H$ is an arbitrary
Hilbert space, and all finite sequences $a_1\,, \ldots \,, a_m\in
A$\,, $b_1\,, \ldots \,, b_n\in B$\,, $c_1\,, \ldots \,, c_m\in
C$\,, $d_1\,, \ldots \,, d_n\in D$\,, where $m, n\in \mathbb{N}$\,,
we have
\begin{equation}\label{eq468}
\left\|\sum\limits_{i=1}^m \sum\limits_{j=1}^n u(a_i,
b_j)\pi(c_i)\rho(d_j)\right\|\leq
2\|u\|_{\text{jcb}}\left\|\sum\limits_{i=1}^m a_i\otimes
c_i\right\|_{A\otimes_{\text{min}} C} \left\|\sum\limits_{j=1}^n
b_j\otimes d_j\right\|_{B\otimes_{\text{min}} D}\,.
\end{equation}
By our main theorem, $\|u\|_{\text{ER}}\leq \|u\|_{\text{jcb}}<
\infty$\,. Let $\xi\,, \eta$ be unit vectors in $H$\,. Let
$u=u_1+u_2$ be the decomposition of $u$ satisfying (\ref{eq456}) and
(\ref{eq457}) as above. Then
\begin{eqnarray}\label{eq469}
&&\!\!\!\!\!\left|\langle\sum\limits_{i=1}^m \sum\limits_{j=1}^n
u(a_i, b_j)\pi(c_i)\rho(d_j)\xi, \eta\rangle\right|\\&&\qquad \leq
\left|\langle \sum\limits_{i=1}^m \sum\limits_{j=1}^n u_1(a_i, b_j)
\pi(c_i)\rho(d_j)\xi, \eta\rangle\right|+ \left|\sum\limits_{i=1}^m
\sum\limits_{j=1}^n \langle u_2(a_i, b_j)\pi(c_i)\rho(d_j)\xi,
\eta\rangle\right|\nonumber\,,
\end{eqnarray}
where $\langle \cdot\,, \cdot \rangle $ denotes the inner product on
$H$\,. By using the GNS construction for the states $f_1$ on $A$ and
$g_1$ on $B$ and inequality (\ref{eq456}), we obtain for any $a\in
A$ and $b\in B$ that
\begin{eqnarray*}
|u_1(a, b)|&\leq &\|u\|_{ER}f_1(aa^*)^{1/2}g_1(b^*b)^{1/2}\\
&=& \|u\|_{ER}\|\pi_{f_1}(a^*)\xi_{f_1}\|\cdot
\|\pi_{g_1}(b)\xi_{g_1}\|\,,
\end{eqnarray*}
where $(H_{f_1}\,, \pi_{f_1}\,, \xi_{f_1})$ is the GNS triple
associated to $(A, f_1)$\,, respectively, $(H_{g_1}\,, \pi_{g_1}\,,
\xi_{g_1})$ is the GNS triple associated to $(B, g_1)$\,. Hence,
there exists $V_1\in {\mathcal B}(H_{g_1}\,, H_{f_1})$ such that
$\|V_1\|\leq \|u\|_{ER}$\,, satisfying
\begin{equation*}
u_1(a, b)=\langle V_1 \pi_{g_1}(b)\xi_{g_1}\,,
\pi_{f_1}(a^*)\xi_{f_1}\rangle\,, \quad a\in A\,, b\in B\,.
\end{equation*}
Therefore, for any $a\in A$ and $b\in B$\,,
\begin{eqnarray}\label{eq470}
&&\!\!\!\!\!\left|\langle \sum\limits_{i=1}^m \sum\limits_{j=1}^n
u_1(a_i, b_j) \pi(c_i)\rho(d_j)\xi, \eta\rangle\right|=
\left|\langle \sum\limits_{i=1}^m \sum\limits_{j=1}^n \langle V_1
\pi_{g_1}(b)\xi_{g_1}\,, \pi_{f_1}(a^*)\xi_{f_1}\rangle
\rho(d_j)\pi(c_i) \xi\,, \eta\rangle\right|\\
&& \qquad =\left|\langle (V_1\otimes 1_H)(\pi_{g_1}\otimes
\rho)(\sum_{j=1}^n b_j\otimes d_j)(\xi_{g_1}\otimes \xi)\,,
(\pi_{f_1}\otimes \pi)(\sum_{i=1}^m a_i^* \otimes
c_i^*)(\xi_{f_1}\otimes\eta)\rangle\right|\nonumber\\
&& \qquad \leq \|u\|_{ER}\left\|\sum_{i=1}^m a_i\otimes
c_i\right\|_{A\otimes_{\text{min}} C} \left\|\sum_{j=1}^n b_j\otimes
d_j\right\|_{B\otimes_{\text{min}} D}\,,\nonumber
\end{eqnarray}
wherein we used the fact that the representations $\pi$ and $\rho$
do commute, and that $\sum_i a_i^*\otimes c_i^*=\left(\sum_i
a_i\otimes c_i\right)^*$\,.

Similarly, by using the GNS construction for the states $f_2$ on $A$
and $g_2$ on $B$ and inequality (\ref{eq457}), we obtain for any
$a\in A$ and $b\in B$ that
\begin{eqnarray*}
|u_2(a, b)|&\leq &\|u\|_{ER}f_2(a^*a)^{1/2}g_2(bb^*)^{1/2}\\
&=& \|u\|_{ER}\|\pi_{f_2}(a)\xi_{f_2}\|\cdot
\|\pi_{g_2}(b^*)\xi_{g_2}\|\,,
\end{eqnarray*}
where $(H_{f_2}\,, \pi_{f_2}\,, \xi_{f_2})$ is the GNS triple
associated to $(A, f_2)$\,, respectively, $(H_{g_2}\,, \pi_{g_2}\,,
\xi_{g_2})$ is the GNS triple associated to $(B, g_2)$\,. Hence,
there exists $V_2\in {\mathcal B}(H_{f_2}\,, H_{g_2})$ such that
$\|V_2\|\leq \|u\|_{ER}$\,, satisfying
\begin{equation*}
u_2(a, b)=\langle V_2 \pi_{f_2}(a)\xi_{f_2}\,,
\pi_{g_2}(b^*)\xi_{g_2}\rangle\,, \quad a\in A\,, b\in B\,.
\end{equation*}
Therefore, for any $a\in A$ and $b\in B$\,,
\begin{eqnarray}\label{eq471}
&&\!\!\!\!\!\left|\langle \sum\limits_{i=1}^m \sum\limits_{j=1}^n
u_2(a_i, b_j) \pi(c_i)\rho(d_j)\xi, \eta\rangle\right|=
\left|\langle \sum\limits_{i=1}^m \sum\limits_{j=1}^n \langle V_2
\pi_{f_2}(a)\xi_{f_2}\,, \pi_{g_2}(b^*)\xi_{g_2}\rangle
\pi(c_i)\rho(d_j) \xi\,, \eta\rangle\right|\\
&& \qquad =\left|\langle (V_2\otimes 1_H)(\pi_{f_2}\otimes
\pi)(\sum_{i=1}^m a_i\otimes c_i)(\xi_{f_2}\otimes \xi)\,,
(\pi_{g_2}\otimes \rho)(\sum_{j=1}^n b_j^* \otimes
d_j^*)(\xi_{g_2}\otimes\eta)\rangle\right|\nonumber\\
&& \qquad \leq \|u\|_{ER}\left\|\sum_{i=1}^m a_i\otimes
c_i\right\|_{A\otimes_{\text{min}} C} \left\|\sum_{j=1}^n b_j\otimes
d_j\right\|_{B\otimes_{\text{min}} D}\,.\nonumber
\end{eqnarray}
The inequality (\ref{eq468}) follows now by (\ref{eq469})\,,
(\ref{eq470}) and (\ref{eq471})\,, since $\|u\|_{\text{ER}}\leq
\|u\|_{\text{jcb}}$\,. The proof is complete.
\end{proof}

Our next proposition gives a complete characterization of those
bilinear forms $u:E\times F\rightarrow \mathbb{C}$ on operator
spaces $E\subseteq A$ and $F\subseteq B$ sitting in C$^*$-algebras
$A$ and $B$\,, for which $\|u\|_{\text{ER}}< \infty$\,.

\begin{prop}\label{prop66664444343432332}
Let $E\subseteq A$ and $F\subseteq B$ be operator spaces sitting in
C$^*$-algebras $A$ and $B$\,, and let $u:E\times F\rightarrow
\mathbb{C}$ be a bounded bilinear map. The following two conditions
are equivalent:
\begin{enumerate}
\item [$(i)$] $\|u\|_{\text{ER}}< \infty$\,.
\item [$(ii)$] There exists a constant $\kappa_4\geq 0$ such that for
all C$^*$-algebras $C$ and $D$\,, all $m, n\in \mathbb{N}$ and all
$a_1\,, \ldots \,,a_m\in E$\,, $b_1\,, \ldots \,,b_n\in F$\,,
$c_1\,, \ldots \,,c_m\in C$\,, $d_1\,, \ldots \,,d_n\in D$\,, we
have
\begin{equation}\label{eq481}
\left\|\sum\limits_{i=1}^m \sum\limits_{j=1}^n u(a_i, b_j)c_i\otimes
d_j \right\|_{C\otimes_{\text{max}} D}\leq \kappa_4
\left\|\sum\limits_{i=1}^m a_i\otimes
c_i\right\|_{E\otimes_{\text{min}} C}\left\|\sum\limits_{j=1}^n
b_j\otimes d_j\right\|_{F\otimes_{\text{min}} D}\,.
\end{equation}
\end{enumerate}
Moreover, if $\kappa_4(u)$ denotes the best constant in $(ii)$\,,
then
\begin{equation*}
\frac12 \|u\|_{\text{ER}}\leq \kappa_4(u)\leq 2\|u\|_{\text{ER}}\,.
\end{equation*}
\end{prop}

\begin{proof}
The implication $(i)\Rightarrow (ii)$ can be obtained from the proof
of Proposition \ref{prop789} with minor modifications. In the case
when $E=A$ and $F=B$ we have by (\ref{eq470}) and (\ref{eq471}) that
\begin{equation}\label{eq77777111111171}
\left\|\sum\limits_{i=1}^m \sum\limits_{j=1}^n u(a_i, b_j)c_i\otimes
d_j \right\|_{C\otimes_{\text{max}} D}\leq 2\|u\|_{\text{ER}}
\left\|\sum\limits_{i=1}^m a_i\otimes
c_i\right\|_{A\otimes_{\text{min}} C}\left\|\sum\limits_{j=1}^n
b_j\otimes d_j\right\|_{B\otimes_{\text{min}} D}\,.
\end{equation}
To extend the proof of (\ref{eq77777111111171}) to the general case
of operator spaces $E\subseteq A$ and $F\subseteq B$\,, the
operators $V_1\in {\mathcal B}(H_{g_1}\,, H_{f_1})$ and $V_2\in
{\mathcal B}(H_{f_2}\,, H_{g_2})$ will instead be operators in
${\mathcal B}(H_{g_1}^0\,, H_{f_1}^0)$ and ${\mathcal
B}(H_{f_2}^0\,, H_{g_2}^0)$\,, respectively, where
\begin{equation*}
H_{f_1}^0:=\overline{\pi_{f_1}(E)^*\xi_{f_1}}\,, \quad
H_{g_1}^0:=\overline{\pi_{g_1}(F)\xi_{g_1}}\,, \quad
H_{f_2}^0:=\overline{\pi_{f_2}(E)\xi_{f_2}}\,, \quad
H_{g_2}^0:=\overline{\pi_{g_2}(F)^*\xi_{g_2}}\,.
\end{equation*}
The rest of the proof of the implication $(i)\Rightarrow (ii)$ can
then be completed as in the Proof of Proposition \ref{prop789}\,. It
also follows that $\kappa_4(u)\leq 2\|u\|_{\text{ER}}$\,.

The converse implication $(ii)\Rightarrow (i)$ can be obtained from
the proof of Theorem 0.3 in \cite{PiSh}\,. For convenience of the
reader, and in order to obtain a better constant, we include below a
slightly modified argument.

Let $u:E\times F\rightarrow \mathbb{C}$ be a bounded bilinear form
satisfying (\ref{eq481}). We will show that $\|u\|_{\text{ER}}\leq
2\kappa_4$\,. By Lemma 2.4 in \cite{PiSh}\,, given a positive
integer $n$ and $\lambda_1\,, \ldots \,, \lambda_n> 0$\,, we can
find two sets $\{x_1\,, \ldots\,, x_n\}$ and $\{y_1\,, \ldots\,,
y_n\}$ of operators on a Hilbert space $H$ with a unit vector
$\Omega$ such that the following properties hold:
\begin{enumerate}
\item [$(a)$] For all $a_1\,, \ldots \,,
a_n\in E$ and all $b_1\,, \ldots \,, b_n\in B$\,,
\begin{eqnarray*}
\left\|\sum_{i=1}^n a_i\otimes x_i\right\|&\leq &\left\|\sum_{i=1}^n
\lambda_i a_i a_i^*\right\|^{1/2}+ \left\|\sum_{i=1}^n
\lambda_i^{-1}a_i^*
a_i\right\|^{1/2}\\
\left\|\sum_{i=1}^n b_i\otimes y_i\right\|&\leq &\left\|\sum_{i=1}^n
\lambda_i b_i b_i^*\right\|^{1/2}+ \left\|\sum_{i=1}^n
\lambda_i^{-1}b_i^* b_i\right\|^{1/2}
\end{eqnarray*}
\item [$(b)$] The von Neumann
algebra $W^*(x_1\,, \ldots\,, x_n)$ generated by $x_1\,, \ldots\,,
x_n$ commutes with the von Neumann algebra $W^*(y_1\,, \ldots\,,
y_n)$ generated by $y_1\,, \ldots\,, y_n$\,.
\item [$(c)$] $\langle x_i y_j \Omega\,,
\Omega\rangle_H=\delta_{ij }$\,, for all $1\leq i, j\leq n$\,.
\end{enumerate}
Let now $n\in \mathbb{N}$ and let $\lambda_1\,, \ldots\,, \lambda_n>
0$ (arbitrarily chosen), be fixed. Then by (\ref{eq481}) we have for
all $a_1\,, \ldots\,, a_n\in E$ and $b_1\,, \ldots \,, b_n\in F$\,,
that
\begin{eqnarray*}
\left|\sum_{i=1}^n u(a_i\,, b_i)\right|&=&\left| \sum\limits_{i,
j=1}^n u(a_i, b_j)\langle
x_iy_j\Omega\,, \Omega\rangle_H\right|\\
&\leq &\left\|\sum\limits_{i, j=1}^n u(a_i,
b_j)x_iy_j\right\|_{{\mathcal B}(H)}\\
&\leq & \left\|\sum\limits_{i, j=1}^n u(a_i, b_j)x_i\otimes
y_j\right\|_{{W^*(x_1\,, \ldots \,, x_n)}\otimes_{\text{max}}
{W^*(y_1\,, \ldots \,, y_n)}}\\
&\leq & \kappa_4 \left\|\sum_{i=1}^n a_i\otimes
x_i\right\|_{E\otimes_{\text{min}} {W^*(x_1\,, \ldots \,, x_n)}}
\left\|\sum_{i=1}^n b_i\otimes
y_i\right\|_{F\otimes_{\text{min}} {W^*(y_1\,, \ldots \,, y_n)}}\\
&\leq & \kappa_4 \left(\left\|\sum_{i=1}^n \lambda_i a_i
a_i^*\right\|^{\frac12}+\left\|\sum_{i=1}^n \lambda_i^{-1}a_i^*
a_i\right\|^{\frac12}\right)\left(\left\|\sum_{i=1}^n \lambda_i
b_i b_i^*\right\|^{\frac12}+\left\|\sum_{i=1}^n \lambda_i^{-1}b_i^* b_i\right\|^{\frac12}\right)\\
&\leq & 2\kappa_4 \left(\left\|\sum_{i=1}^n \lambda_i a_i
a_i^*\right\|+\left\|\sum_{i=1}^n \lambda_i^{-1}a_i^*
a_i\right\|\right)^{\frac12}\left(\left\|\sum_{i=1}^n \lambda_i b_i
b_i^*\right\|+\left\|\sum_{i=1}^n \lambda_i^{-1}b_i^*
b_i\right\|\right)^{\frac12}\,,
\end{eqnarray*}
where we have used the well-known inequality
\[ \sqrt{\alpha}+\sqrt{\beta}\leq \sqrt{2}\sqrt{\alpha+\beta}\,, \quad \alpha\,,
\beta\geq 0\,. \] Since $2\sqrt{\alpha \beta}\leq \alpha+ \beta$ for
all $\alpha\,, \beta\geq 0$\,, it follows that
\begin{equation}\label{eq491}
\left|\sum_{i=1}^n u(a_i\,, b_i)\right|\leq \kappa_4
\left(\left\|\sum_{i=1}^n \lambda_i a_i^* a_i\right\|+
\left\|\sum_{i=1}^n \lambda_i^{-1} a_i a_i^*\right\|+
\left\|\sum_{i=1}^n \lambda_i b_i b_i^*\right\|+\left\|\sum_{i=1}^n
\lambda_i^{-1} b_i^* b_i\right\|\right)\,.
\end{equation}
Using a Pietsch separation argument similar to the one given in the
proof of Lemma 3.4 in \cite{Ha6}, we infer the existence of states
$f_1\,, f_2$ on $A$ and $g_1\,, g_2$ on $B$ such that for all $a\in
E$\,, $b\in F$ and $\lambda> 0$\,,
\begin{equation*}
|u(a, b)|\leq \kappa_4 \left(\lambda f_1(aa^*)+\lambda^{-1}
f_2(a^*a)+ \lambda g_2(bb^*)+\lambda^{-1} g_1(b^*b)\right)\,.
\end{equation*}
Replacing now $a$ by $t^{1/2}a$ and $b$ by $t^{{-1}/2}b$\,, where
$t> 0$\,, it follows that for all $a\in E$\,, $b\in F$\,, $t> 0$ and
$\lambda> 0$\,,
\begin{equation*}
|u(a, b)|\leq \kappa_4 \left(t\lambda f_1(aa^*)+t\lambda^{-1}
f_2(a^*a)+\frac1{t} \lambda g_2(bb^*)+\frac1{t} \lambda^{-1}
g_1(b^*b)\right)\,.
\end{equation*}
By taking the infimum over all $t> 0$\,, we deduce by
(\ref{eq22222452}) that for all $a\in E$\,, $b\in F$ and all
$\lambda> 0$\,,
\begin{equation*}
|u(a, b)|\leq 2\kappa_4(\lambda f_1(aa^*)+\lambda^{-1}
f_2(a^*a))^{1/2}(\lambda g_2(bb^*)+\lambda^{-1} g_1(b^* b))^{1/2}\,.
\end{equation*}
Therefore, for all $a\in E$\,, $b\in F$ and $\lambda> 0$\,,
\begin{equation*}
|u(a, b)|^2\leq
(2\kappa_4)^2(f_1(aa^*)g_1(b^*b)+f_2(a^*a)g_2(bb^*)+\lambda^2
f_1(aa^*)g_2(bb^*)+\lambda^{-2} f_2(a^*a)g_1(b^*b))\,.
\end{equation*}
By taking infimum over $\lambda> 0$\,, a further application of
(\ref{eq22222452}) shows that for all $a\in E$ and $b\in F$\,,
\begin{eqnarray*}
|u(a, b)|^2&\leq &
(2\kappa_4)^2\left(f_1(aa^*)g_1(b^*b)+f_2(a^*a)g_2(bb^*)+2f_1(aa^*)^{1/2}g_1(b^*b)^{1/2}f_2(a^*a)^{1/2}g_2(bb^*)^{1/2}\right)\\
&=&
(2\kappa_4)^2\left(f_1(aa^*)^{1/2}g_1(b^*b)^{1/2}+f_2(a^*a)^{1/2}g_2(bb^*)^{1/2}\right)^2\,.
\end{eqnarray*}
This implies that $\|u\|_{ER}\leq 2\kappa_4$\,, which completes the
proof of the implication $(ii)\Rightarrow (i)$ and it also proves
the inequality $\kappa_4\geq \frac12 \|u\|_{\text{ER}}$\,.
\end{proof}

\section*{Acknowledgements} This paper was completed during the
authors stay at the Fields Institute in the Fall of 2007, while
attending the Thematic Program on Operator Algebras. We would like
to thank the Fields Institute and the organizers of the program for
their support and warm hospitality.

\thanks{}


\vspace*{0.3cm}
\begin{thebibliography}{456}
\bibitem{Bl}{\sc D. Blecher}, {\em Generalizing Grothendieck's
program}, Function spaces, Edited by K. Jarosz, Lect. Notes in Pure
and Applied Math., Vol. 136, Marcel Dekker, 1992.
\bibitem{ChS}{\sc E. Christensen, A. Sinclair}, {\em A survey of
completely bounded operators}, Bull. London Math. Soc. {\bf 21}
(1989), 417-448.
\bibitem{Co4}{\sc A. Connes}, {\em Une classification des facteurs
de type III}, Ann. Scient. \'Ec. Norm. Sup., $4^{\text{e}}$ s\'erie,
tome 6 (1973), 133-252.
\bibitem{Con}{\sc A. Connes}, {\em Classification of injective factors. Cases II$_1$\,, II$_\infty$\,, III$_\lambda$\,, $\lambda\ne 1\,.$}, Ann. Math. {\bf 104} (1976), 73-115.
%\bibitem{Co7}{\sc A. Connes}, {\em Factors of type III$_1$, property $L_\lambda'$ and closure of inner automorphisms}, J. Op. Theory {\bf 14} (1985), no. 1, 189--211.
\bibitem{ER2}{\sc E. Effros and  Z.-J. Ruan}, {\em A new approach to
operator spaces}, Canad. Math. Bull. Vol. {34} (3) (1991), 329-337.
\bibitem{ER3}{\sc E. Effros and  Z.-J. Ruan}, {\em Self-duality for
the Haagerup tensor product and Hilbert space factorization}, J.
Funct. Analysis {\bf 100} (1991), 257-284.
\bibitem{ER}{\sc E. Effros and  Z.-J. Ruan}, {\em Operator Spaces}, London Math. Soc. Monographs New Series {\bf 23}, Oxford University Press, 2000.
\bibitem{Gro}{\sc A. Grothendieck}, {\em Resum\'e de la th\'eorie
m\'etrique des produits tensorielles topologiques}, Bol. Soc. Mat.
Sao Paolo {\bf 8} (1956), 1-79.
%\bibitem{Ha7}{\sc U. Haagerup}, {\em Solution of the similarity
%problem for cyclic representations of C$^*$-algebras}, Ann. of Math.
%{\bf 118} (1983), 215-240.
\bibitem{Ha6}{\sc U. Haagerup}, {\em The Grothendieck inequality for
bilinear forms on C$^*$-algebras}, Adv. Math., Vol. 56, No. 2
(1985), 93-116.
%\bibitem{Ha5}{\sc U. Haagerup}, {\em A new proof of the equivalence
%of injectivity and hyperfiniteness for factors on a separable
%Hilbert space}, J. Funct. Analysis, Vol. 62, No. 2 (1985), 160-201.
%\bibitem{Ha2}{\sc U. Haagerup}, {\em Connes' bicentralizer problem and uniqueness of the injective factor of type $III_1$}, Acta Math. {\bf 158} (1987), no.1-2, 95-148.
\bibitem{Ha1}{\sc U. Haagerup}, {\em The injective factors of type
III$_\lambda$\,, $0< \lambda< 1$}, Pacific J. Math., Vol. 137, No. 2
(1989), 265-310.
%\bibitem{HM}{\sc U. Haagerup, M. Musat}, {\em On the best constants
%in noncommutative Khintchine-type inequalities}, J. Funct. Analysis,
%to appear.
%\bibitem{Ju2}{\sc M. Junge}, {\em Embedding of the operator space
%$OH$ and the logarithmic "little Grothendieck inequality"}, Invent.
%Math. {\bf 161} (2005), 389-406.
%\bibitem{Ju}{\sc M. Junge}, {\em Operator spaces and Araki-Woods factors-A quantum probabilistic approach}, Int. Math. Res. Pap. 2006, Art. ID 76978, 87 pp.
%\bibitem{KR}{\sc R. V. Kadison and J. R. Ringrose}, {\em Fundamentals of the Theory of Operator Algebras I, II}, Academic Press, 1986.
%\bibitem{Ki}{\sc E. Kirchberg}, {\em On nonsemisplit extensions,
%tensor products and exactness of group $C^*$-algebras}, Invent. Math
%{\bf 112} (1993), no. 3, 449-489.
\bibitem{LT}{\sc J. Lindenstrauss and L. Tzafriri}, {\em Classical
Banach Spaces, Sequence Spaces}, Ergebnisse, Vol. {92},
Springer-Verlag, 1992.
%\bibitem{Pau}{\sc V. I. Paulsen}, {\em Completely Bounded Maps and Dilations}, Pitman Res. Notes, Longman Sci. Tech., London, 1986.
\bibitem{Pi7}{\sc G. Pisier}, {\em Grothendieck's theorem for
non-commutative C$^*$-algebras with an appendix on Grothendieck's
constant}, J. Funct. Analysis {\bf 29} (1978), 397-415.
\bibitem{Pi6}{\sc G. Pisier}, {\em The Operator Hilbert Space $OH$,
Complex Interpolation and Tensor Norms}, Mem. Amer. Math. Soc.
Number {\bf 585}, Vol. 122, Providence, RI, 1996.
\bibitem{Pi1}{\sc G. Pisier}, {\em An Introduction to the Theory of Operator Spaces}, London Math. Soc. Lect. Notes Series {\bf 294}, Cambridge University Press, Cambridge 2003.
%\bibitem{Pi2}{\sc G. Pisier}, {\em Completely bounded maps into certain Hilbertian operator spaces}, Int. Math. Res. Not. (2004), no.74, 3983-4018.
%\bibitem{Pi3}{\sc G. Pisier}, {\em The operator Hilbert space OH and type III von Neumann algebras}, Bull. London Math. Soc. {\bf 36} (2004), no.4, 455-459.
\bibitem{PiSh}{\sc G. Pisier and D. Shlyakhtenko}, {Grothendieck's
theorem for operator spaces}, Invent. Math. {\bf 150} (2002),
185-217.
%\bibitem{Sak}{\sc S. Sakai}, {\em $C^*$-algebras and
%$W^*$-algebras}, Springer Verlag Berlin, Heidelberg, 1998.
\bibitem{Ta1}{\sc M. Takesaki}: {\em The structure of von Neumann algebras with a homogeneous periodic state}, Acta
Math. {\bf 131} (1973), 79-121.
\bibitem{Tk}{\sc M. Takesaki}, {\em Theory of Operator Algebras II, III}, Springer-Verlag, New-Yord, 1979.
\bibitem{Xu}{\sc Q. Xu}, {\em Operator space Grothendieck inequalities for noncommutative $L_p$-spaces}, Duke Math. J. {\bf
131} (2006), 525-574.

\end{thebibliography}
\end{document}